\documentclass[12pt]{article}
\usepackage{amsfonts}
\usepackage{amsmath}
\usepackage{amssymb}
\usepackage{amscd}
\usepackage{eepic}
\usepackage{epic}
\usepackage[francais]{babel}
\usepackage[latin1]{inputenc}
\usepackage{color}
\usepackage{graphics}
\begin{document}
\newtheorem{Theoreme}{Th\'eor\`eme}[section]
\newtheorem{Theorem}{Theorem}[section]
\newtheorem{Th}{Th\'eor\`eme}[section]
\newtheorem{De}[Th]{D\'efinition}
\newtheorem{Pro}[Th]{Proposition}
\newtheorem{Lemma}[Theorem]{Lemma}
\newtheorem{Proposition}[Theoreme]{Proposition}
\newtheorem{Lemme}[Theoreme]{Lemme}
\newtheorem{Corollaire}[Theoreme]{Corollaire}
\newtheorem{Consequence}[Theoreme]{Cons\'equence}
\newtheorem{Remarque1}[Theoreme]{Remarque}
\newtheorem{Convention}[Theoreme]{{\sc Convention}}
\newtheorem{PP}[Theoreme]{Propri\'et\'es}
\newtheorem{Conclusion}[Theoreme]{Conclusion}
\newtheorem{Ex}[Theoreme]{Exemple}
\newtheorem{Definition}[Theoreme]{D\'efinition}
\newtheorem{Remark1}[Theorem]{Remark}
\newtheorem{Not}[Theoreme]{Notation}
\newtheorem{Nota}[Theorem]{Notation}
\newtheorem{Propo}[Theorem]{Proposition}
\newtheorem{exercice1}[Th]{Lemme-Confi� au lecteur}
\newtheorem{Corollary}[Theorem]{Corollary}
\newtheorem{PPtes}[Th]{Propri\'et\'es}
\newtheorem{Defi}[Theorem]{Definition}
\newtheorem{Example1}[Theorem]{Example}
\newenvironment{Proof}{\medbreak{\noindent\bf Proof }}{~{\hfill $\bullet$\bigbreak}}

\newenvironment{Demonstration}{\medbreak{\noindent\bf D\'emonstration
 }}{~{\hskip 3pt$\bullet$\bigbreak}} 

\newenvironment{Remarque}{\begin{Remarque1}\em}{\end{Remarque1}} 
\newenvironment{Remark}{\begin{Remark1}\em}{\end{Remark1}}
\newenvironment{Exemple}{\begin{Ex}\em}{~{\hskip
3pt$\bullet$}\end{Ex}} 
\newenvironment{exercice}{\begin{exercice1}\em}{\end{exercice1}}
\newenvironment{Notation}{\begin{Not}\em}{\end{Not}}
\newenvironment{Notation1}{\begin{Nota}\em}{\end{Nota}}
\newenvironment{Example}{\begin{Example1}\em}{~{\hskip
3pt$\bullet$}\end{Example1}} 

\newenvironment{Remarques}{\begin{Remarque1}\em \ \\* }{\end{Remarque1}}

\renewcommand{\Re}{{\cal R}}
\newcommand{\Dim}{{\rm Dim\,}}
\renewcommand{\Im}{{\cal F}}
\newcommand{\finpreuve}{~{\hskip 3pt$\bullet$\bigbreak}}
\newcommand{\hp}{\hskip 3pt}
\newcommand{\hph}{\hskip 8pt}
\newcommand{\hphh}{\hskip 15pt}
\newcommand{\vp}{\vskip 3pt}
\newcommand{\vpv}{\vskip 15pt}
\newcommand{\IP}{{\mathbb{IP}}}
\newcommand{\rd}{{\mathbb R}^2}
\newcommand{\R}{{\mathbb R}}

\newcommand{\Hyper}{{\mathbb H}}
\newcommand{\Int}{{\mathbb I}}
\newcommand{\Boule}{{\mathbb B}(0,1)}
\newcommand{\Cantor}{{\mathbb K}}
\newcommand{\Sp}{{\mathbb S}}
\newcommand{\K}{{\mathbb K}}
\newcommand{\B}{{\mathbb B}}
\newcommand{\Ho}{{\mathbb H}}
\newcommand{\Nat}{{\mathbb N}}
\newcommand{\N}{{\mathbb N}}
\newcommand{\Proba}{{\mathbb P}}
\newcommand{\Esp}{{\mathbb E}}
\newcommand{\Complex}{{\mathbb C}}
\newcommand{\Ha}{{\cal H}}
\newcommand{\Harm}{{\bold H}}
\newcommand{\Lcal}{{\cal L}}
\newcommand{\ds}{\displaystyle}
\newcommand{\un}{\bold 1}
\newcommand{\Cone}{C(x,r,\epsilon ,\Phi)}
\newcommand{\Cn}{C(x,2^{-n},\epsilon ,\Phi )}
\newcommand{\Tranche}{W(x,r,\epsilon,\Phi)}
\newcommand{\Wn}{W(x,2^{-n},\epsilon,\Phi)}
\newcommand{\WFn}{W(x,2^{-n},\epsilon,\Phi)\cap F}
\newcommand{\ovec}{\overrightarrow}
\newcommand{\red}{{\bold R}}
\newcommand{\dimH}{\dim_{\Ha}}
\newcommand{\diam}{\mbox{diam}}
\newcommand{\diamit}{\mbox{\em diam}}
\newcommand{\para}{\vskip 2mm}
\newcommand{\cod}{\stackrel{\mbox{\tiny cod}}{\sim}}
\newcommand{\cardit}{\mbox{\em card}}
\newcommand{\card}{\mbox{card}}
\newcommand{\Sphere}{{\mathbb S}_d}
\newcommand{\distit}{\mbox{\em dist}}
\newcommand{\Tri}{{\cal P}}
\newcommand{\LL}{{\mathcal L}}
\newcommand{\infess}{\mbox{inf\,ess}}
\newcommand{\supess}{\mbox{sup\,ess}}

\definecolor{darkblue}{rgb}{0,0,.5}
\def\u{\underline }
\def\o{\overline}
\def\h{\hskip 3pt}
\def\hh{\hskip 8pt}
\def\hhh{\hskip 15pt}
\def\v{\vskip 8pt}
\def\vv{\vskip 15pt}
\font\courrier=cmr12
\font\grand=cmbxti10
%\font\larger=cmbxti10
\font\large=cmbx12
\font\largeplus=cmr17
\font\small=cmbx8
\font\nor=cmbxti10
\font\smaller=cmr8
\font\smallo=cmbxti10
\openup 0.3mm
%\setbox1\vbox{\hspace{.1cm}{\Large \bf\color{darkblue}Athanasios BATAKIS}%
%
%{\small\em \noindent \hskip 10mm  D\'epartement de Math\'ematiques
%\vskip-2mm
%\noindent \hskip 20mm  Universit� d'Orl�ans
%\vskip-2mm
%\noindent \hskip 18mm  45067 Orl\'eans cedex 2
%\vskip-2mm
%\noindent \hskip 21mm T\'el : 02.38.41.73.15}
%
%\noindent e-mail : Athanasios.Batakis@univ-orleans.fr}
%\ht1=12mm \dp1=10mm \wd1=0mm
%
%\setbox3\hbox {\box1  \hskip 115mm}
%\
%\vskip -1cm \hskip -2cm \box3

\title{Multifractal Analysis of inhomogeneous Bernoulli products}
\author{Athanasios BATAKIS and Beno\^\i t TESTUD}
\date{}
\maketitle
{\bf Abstract} We are interested to the multifractal analysis of  inhomogeneous Bernoulli products  which are also known as coin tossing measures. We give conditions ensuring the validity of the  multifractal formalism for such measures. On another hand, we show that these measures can have a dense set of phase transitions. 
 

\section{Introduction}

Let us consider the dyadic tree $\mathbb T$ (even though all the results in this paper can be easily generalised to any $\ell$-adic structure, $\ell\in\N$), let $\Sigma=\{0,1\}^{\N}$ be its limit (Cantor) set and denote by $(\Im_n)_{n\in\N}$ the associated filtration with the usual $0-1$ encoding. 

For $\epsilon_1,...,\epsilon_n\in\{0,1\}$ we denote by $I_{\epsilon_1...\epsilon_n}$ the cylinder of the $n$th generation defined by 
$I_{\epsilon_1...\epsilon_n}=\{x=(i_1,...,i_n,i_{n+1},...)\in \Sigma\; ;\;i_1=\epsilon_1,...,i_n=\epsilon_n\}$. For every $x\in \Sigma$, $I_n(x)$ stands for   the cylinder  of  $\mathcal{F}_n$ containing~$x$.

If $(p_n)_n$ is a sequence of weights, $p_n\in (0,1)$,  we are interested in Borel measures $\mu$ on $\Sigma$  defined 
in the following way 
\begin{equation}\label{inv}
\mu(I_{\epsilon_1...\epsilon_n})=\prod_{j=1}^n p_j^{1-\epsilon_j}(1-p_j)^{\epsilon_j}.
\end{equation}

A measure of this form will be referred to as an {\em inhomogeneous Bernoulli product}.
The aim of this paper is to study multifractal properties of such measures.

The particular case where the sequence $(p_n)$ is constant is well-known and provides an example of measure satisfying the multifractal formalism (see e.g \cite{fal}). In the general case,    Bisbas    \cite{Bis} gave a sufficient condition on the sequence $(p_n)$  ensuring that $\mu$ is a multifractal measure (i.e. the level sets are not empty). However, the work of Bisbas does not provide the dimension of the level sets $E_\alpha$ associated to  the measure~$\mu$. 

Let us give a brief description of multifractal formalism. 
For  a probability measure~$m$ on $\Sigma$,  we  define the {\em local dimension} (also called   H\"older exponent) of $m$ at $x\in \Sigma$ by $$\displaystyle \alpha(x)=\liminf_{n\rightarrow +\infty} \alpha_n(x)=\liminf_{n\rightarrow +\infty} -\frac{\log m(I_n(x))}{n\log 2}.$$   The  aim of  multifractal analysis is  to find the Hausdorff dimension, $\dim(E_\alpha)$,  of  the  level  set  $ E_\alpha=\left\{x : \alpha(x)=\alpha\right\}$ for  $\alpha>0$.   The function  $f(\alpha)=\dim(E_\alpha)$ is called the {\em singularity spectrum} (or multifractal spectrum) of $m$ and we say that $m$ is a {\em multifractal measure} when $f(\alpha)>0$ for several $\alpha{}'s$.

The concepts underlying the  multifractal decomposition of a measure  go back to an early  paper of    Mandelbrot \cite{mandel}.   In the 80's multifractal measures were used    by physicists to study various  models arising from natural phenomena.   In    fully developped turbulence they  were  used by Frisch and Parisi \cite{fp} to investigate the intermittent behaviour in the regions of high vorticity. In dynamical system theory they were used by Benzi et al.~\cite{ben}  to measure how often a given region of the attractor is visited.  In diffusion-limited aggregation (DLA) they were  used by Meakin et al.~\cite{mea}  to describe the probability of a random walk landing to the neighborhood of a given site on the aggregate. 

In order to determine the function $f(\alpha)$,       Hentschel and  Procaccia  \cite{hp}  used  ideas based on  Renyi entropies \cite{re} to  introduce the generalized dimensions  $D_q$ defined by  \[D_q= \lim_{n\rightarrow +\infty} \frac{1}{q-1}  \frac{\log\left(\sum_{I\in\mathcal{F}_n} m(I)^q\right)}{n \log 2},\] (see also \cite{grp,gr}). From a physical and heuristical point of view, Halsey et al. \cite{has} showed that the singularity spectrum $f(\alpha)$ and the generalized dimensions $D_q$ can be derived from each other. The Legendre transform  turned out to be a useful tool linking  $f(\alpha)$ and $D_q$.  More precisely,  it was suggested that   
\begin{eqnarray}\label{fm} f(\alpha)=\dim(E_{\alpha})=\tau^*(\alpha)=\inf(\alpha q+\tau(q),\ q\in \R),\end{eqnarray}  where \[\tau (q) = \limsup_{n\rightarrow +\infty}\,\tau _n (q)\quad\mbox{ with }\quad \tau _n (q)=\frac{1}{n\log 2}\log\left(\sum_{I\in\mathcal{F}_n} m(I)^q\right).\]   (The sum  runs over the cylinders  $I$  such that $m(I)\not=0$.)   The  function $\tau(q)$ is  called the  $L^q$-spectrum of $m$   and if the limit exists $\tau(q)=(q-1)D_q$. 

Relation \eqref{fm} is called the multifractal formalism and in many aspects it is analogous to the well-known thermodynamic formalism developed by Bowen \cite{bo} and Ruelle \cite{ru}. In general, the main problem is to obtain the minoration $\dim(E_{\alpha})\ge \tau^*(\alpha)$. 

For number of measures, this formalism    can  be verified rigorously. In particular, if the sequence $(p_n)$ is constant or periodic, the measure $\mu$ given by  \eqref{inv}  satisfies  the multifractal formalism (e.g. \cite{fal}). It is also the case  for  invariant measures in some  dynamical systems (e.g \cite{col, Fan,rand}), for  self-similars measures under separation conditions (e.g \cite{cm, Feng,  lng, ol, ri, yyl})  and for quasiindependent measures(e.g \cite{BMP,heurt,test}).

     Despite all these  investigations mentioned, the exact range of the validity of the multifractal formalism is still not known. Olsen \cite{ol} give a rigorous approach  of multifractal formalism in a general context.   This work and the paper of  Brown, Michon
and Peyri\`ere \cite{BMP} enlighten the link between the  minoration $\dim(E_{\alpha})\ge \tau^*(\alpha)$ and the existence of auxiliary measures $m_q$ (the so-called    {\em Gibbs measure} \cite{mi})  satisfying \begin{eqnarray*}\forall n,\,\,\forall I\in \mathcal{F}_n, \quad
  \frac1C  m(I)^q 2^{-n\tau(q)}\le  m_q(I)\le C m(I)^q
  2^{-n\tau(q)},\end{eqnarray*} where the constant  $C>0$ is
independent of $n$ and $I$. In fact, it is shown in \cite{bn, bnbh} that the existence of a measure $m_q$ satisfying  \begin{eqnarray*} m_q(I)\le C m(I)^q
  2^{-n\tau(q)},\end{eqnarray*} is sufficient to obtain the minoration  $\dim(E_{\alpha})\ge \tau^*(\alpha)$ for $\alpha=-\tau'(q)$.     
In this situation, the values of $\alpha$ for which the multifractal formalism may fail lie in intervals $(-\tau'(q^+), \,-\tau'(q^-))$ where $q$ is a point of non differentiability of $\tau$ ($\tau'(q^+)$ and $\tau'(q^-)$ stands for the right and the left derivatives respectively). Such a point $q$ will be called a {\em phase transition}. 
 
If  the weights $p_n$ are not all  the same,  the measure $\mu$ defined by \eqref{inv} is in general no shift-invariant and we cannot apply  classical tools  of ergodic theory, as Shannon-McMillan theorem (e.g \cite{Bil}), to get a lower bound of $\dim(E_\alpha)$ and the differentiability of the function $\tau$.

 Let us introduce the other following level sets defined by
$$\underline E_{\alpha}=\left\{x\; ;\;
\alpha(x)\le \alpha\right\},  \;   \overline F_{\alpha}=\left\{x\; ;\;
\limsup_{n\to\infty}\alpha_n(x)\ge \alpha\right\},$$ and $$F_{\alpha}=\left\{x\; ;\;
\limsup_{n\to\infty}\alpha_n(x)= \alpha\right\}.$$
 
We can now state our main results. In section 2, we prove the following.  
\begin{Theorem}\label{1}
Let $\mu$ be an inhomogeneous Bernoulli product on $\Sigma$ and $q\in\R$.
 We have
$$
\liminf_{n\to\infty}-q\tau_{n}'(q)+\tau_{n}(q)\le\dim\left(\underline E_{-\tau'(q^-)}\cap \overline F_{-\tau'(q^+)}\right)\le
\sup\left\{\tau^*(-\tau'(q^+)), \tau^*(-\tau'(q^-))\right\}.
$$  \end{Theorem}

The proof of the lower bound relies on the construction of a special  inhomogeneous  Bernoulli product which has the dimension of the studied level set.

In section 3 we study  the case  $\alpha=-\tau'(q)$. The existence of $\tau'(q)$ is not sufficient to ensure  the validity of  multifractal formalism  for such values of $\alpha's.$     However, we prove that the multifractal formalism holds if the sequence $(\tau_{n}(q))$ converges. More precisely, we have
\begin{Theorem}\label{4}
Suppose that the sequence $(\tau_{n}(q))$ converges at a point $q\in\R$. If $\tau'(q)$ exists and if $\alpha=-\tau'(q)$, we have 
\begin{equation}
\dim \left(E_{\alpha}\cap F_{\alpha}\right)=\tau^*(\alpha)=\alpha q+\tau(q).
\end{equation}
\end{Theorem}

We easily deduce the following
\begin{Corollary}
 If $p_n$ tend to $p$, as $n\to\infty$, then $\tau_{\mu}=\tau(p,.)$ and the mesure $\mu$ satisfies the multifractal formalism. Nevertheless, the measure $\mu$ can be singular with respect to the (homogeneous) Bernoulli measure associated to $p$.
\end{Corollary}

Theorem \ref{4} leads us to study the differentiability  of the $L^q$-spectrum $\tau(q)$.   In  section 4, we will see that the $L^q$-spectrum of an inhomogeneous Bernoulli product    may be a very irregular function. In particular,  

\begin{Theorem}\label{3}
There exist inhomogeneous Bernoulli products   presenting a dense set of  phase transitions on $(1,+\infty)$.
\end{Theorem}

The are several examples of measures presenting phase transitions (see for instance \cite{testud} and the references therein). The example we propose in this work differs from previous ones at three points: first the phase transitions are situated at points $q>1$ and not at negative ones, where constructions are easier to carry out. Secondly, the set of transitions is dense in $[1,\infty)$, that means as \og bad\fg\  as can be. And finally, the measure presenting this pathologie is just a Bernoulli product! Let us also point out that with some minor modifications our method can also apply to create a dense set of phase transitions within $(0,1)$.

%and if $\alpha=-\tau'(q^-)$
%\begin{equation}
%\liminf_{n\to\infty}-q\tau_{\mu,n}'+\tau_{\mu,n}(q)\le\dim E_{\alpha}\le\tau^*(\alpha).
%\end{equation}

\section{Proof of Theorem \ref{1}}
We begin by a preliminary result.
\begin{Lemma}\label{secderbdd}
If $\mu$ is an inhomogeneous Bernoulli product, then the functions $(\tau_{\mu,n}'')$
 are locally uniformly bounded on $(0,+\infty)$.
 \end{Lemma} 
 \begin{Proof}
 We denote by $\beta(p_i)$ the homogeneous  Bernoulli  measure of parameter $p_i$ and by $\tau(p_i,q)$ it's $\tau$ function, 
 $\tau(p_i,q)=\log_2(p_i^q+(1-p_i)^q)$. Using the fact that $\mu$ is the product of $\beta(p_i)$ we easily obtain
 $$\tau_{\mu,n}(q)=\frac1n\sum_{i=1}^n\tau(p_i,q).$$ 

 It is therefore sufficient to show that, for any $q_0>0$, there exists a constant $C=C(q_0)$ such that 
for all $p\in(0,1)$ and   all $q>q_0$, $\displaystyle\frac{\partial^2\tau(p,q)}{\partial q^2}\le C$.
We have 

\begin{eqnarray*}
\frac{\partial^2\tau(p,q)}{\partial q^2}&=&
\frac{p^q(\log_2 p)^2+(1-p)^q(\log_2(1-p))^2}{p^q+(1-p)^q}-\frac{\left(p^q\log_2 p+(1-p)^q\log_2(1-p)\right)^2}{(p^q+(1-p)^q)^2}\\
&=&\frac{p^q(1-p)^q\left((\log_2 p)^2+(\log_2(1-p))^2-2\log_2 p\log_2(1-p)\right)}{(p^q+(1-p)^q)^2}\\
&=&\frac{p^q(1-p)^q\left(\log_2\frac{p}{1-p}\right)^2}{(p^q+(1-p)^q)^2}\le [4p(1-p)]^q(\log_2\frac{p}{1-p})^2\\&\le&[4p(1-p)]^{q_0}(\log_2\frac{p}{1-p})^2,
\end{eqnarray*}
which is uniformly bounded on $p\in(0,1)$ and the proof is complete.
\end{Proof}

Lemma \ref{secderbdd} allows us to give estimates for the lower and the upper Hausdorff dimension of the measure $\mu$. They are respectively defined by 
$$\dim_*(\mu)=\inf\{\dim(E),\,\,\, \mu(E)>0\} \,;\,\,\dim^*(\mu)=\inf\{\dim(E),\,\,\, \mu(E)=1\}.$$  We say that $\mu$ is exact if $\dim_*(\mu)=\dim^*(\mu)$ and we note $\dim(\mu)$ the common value. In the same way,  we can define  the lower and the upper Packing  dimension $\Dim$ of the measure $\mu$.  It is well known that there exist some relations between these quantities and the derivatives of the function $\tau_\mu(q)$ at $q=1$. More precisely, it is proved in \cite{Fan,heurt} that 
$$-\tau_\mu'(1+)\le\dim_*(\mu)\le  h_*(\mu)\le h^*(\mu)\le \Dim^*(\mu)\le -\tau_\mu'(1-),$$ where $ h_*(\mu)$ and $ h^*(\mu)$ stand for the lower and the upper entropy of the measure $\mu$, defined as 
$$ h_*(\mu)=\liminf -\frac{1}{n\log 2} \sum_{I\in\Im_n} \mu(I)\log \mu(I)= \liminf -\tau_{\mu_n}'(1)$$ and  $$ h^*(\mu)=\limsup -\frac{1}{n\log 2} \sum_{I\in\Im_n} \mu(I)\log \mu(I)= \limsup - \tau_{\mu_n}'(1).$$  

By Lemma \ref{secderbdd}, we deduce (see \cite{BH,heurt}) the following remark.
 \begin{Remark}\label{Heurt}
 If $\mu$ is an inhomogeneous Bernoulli product then 
$$\displaystyle\dim\mu =-\tau_{\mu}'(1^+)=h_*(\mu)=\liminf_{n\to\infty}-\tau_{\mu_n}'(1)$$ and 
$$\displaystyle\Dim\mu=-\tau_{\mu}'(1^-)=h^*(\mu)=\limsup_{n\to\infty}-\tau_{\mu_n}'(1).$$
 \end{Remark}

Fix $q\in \R$. To prove Theorem \ref{1}, we construct an auxiliary  measure $\nu$ supported by the set $\underline E_{-\tau'(q^-)}\cap \overline F_{-\tau'(q^+)}$. More precisely, we consider a sequence of measures  $\nu_n$ satisfying 
\begin{eqnarray}\label{can} \forall I\in\Im_n,\quad  \nu_n(I)=\frac{\mu(I)^q}{\sum_{I\in\Im_n}\mu(I)^q}=\mu(I)^q|I|^{\tau_{\mu,n}(q)}.\end{eqnarray}   ($\vert I \vert = 2^{-n}$ stands for  the diameter of $I$).  The following lemma implies that  the sequence $(\nu_n)$ converges in the weak$^*$ sense to a probability  measure  $\nu$ which 
is by construction  an inhomogeneous Bernoulli product. 
 \begin{Lemma}\label{nuthing}
Let $n\in\N$ and $I\in \Im_n$. If  $\mu$ is an inhomogeneous Bernoulli product, we have 
$\nu_n(I) =\nu_{n+1}(I)$. \end{Lemma}
\begin{Proof}
Take $n>0$ and $I\in\Im_n$. We can compute
$$\nu_{n+1}(I)=\frac{\sum_{J\in\Im_1}\mu(IJ)^q}{\sum_{I\in\Im_n}\sum_{J \in\Im_1}\mu(IJ)^q}=
\frac{\mu(I)^q(p_{n+1}^q+(1-p_{n+1})^q)}{\sum_{I\in\Im_n}(p_{n+1}^q+(1-p_{n+1})^q)\mu(I)^q}$$
and therefore
$\nu_{n+1}(I)=\nu_n(I)$ for all $I\in \Im_n$.
\end{Proof}
By remark \ref{Heurt}, we then deduce that the Hausdorff and the Packing dimension of $\nu$ are given by an entropy formula. In other terms, we have $$\displaystyle\dim\nu=\liminf_{n\to\infty}-\tau_{\nu,n}'(1)=h_*(\nu)$$
and
$$\displaystyle\Dim\nu=\limsup_{n\to\infty}-\tau_{\nu,n}'(1)=h^*(\nu).$$

Now we can  prove Theorem \ref{1}.

 \begin{Proof}{\bf of Theorem \ref{1}}
The upper bound is a well known fact of multifractal formalism (see for instance \cite{BMP}). 
In fact we have 
\begin{enumerate}
\item If $\alpha\le -\tau'(0^+)$ then 
$\dim E_{\alpha}\le\dim \underline E_{\alpha}\le\tau^*(\alpha).$
\item If $\alpha\ge -\tau'(0^-)$ then 
$\dim F_{\alpha}\le\dim \overline F_{\alpha}\le\tau^*(\alpha).$
\item $-\tau'(0^+)\le \alpha\le -\tau'(0^-)$ then $\tau^*(\alpha)=\tau(0)=1$ and the upper bound follows.
\end{enumerate}

 Relation \eqref{can} easily gives
$\tau_{\nu,n}(s)=\tau_{\mu,n}(qs)-s\tau_{\mu,n}(q)$.
From  remark \ref{Heurt}, using  the inhomogeneous Bernoulli property of $\mu$ and $\nu$,    we deduce   that 
%$\tau_{\nu}'(1^-)=\limsup\tau_{\nu,n}'(1)=\limsup\left(q\tau_{\mu,n}'(q)-\tau_{\mu,n}(q)\right)$ and similarly
$$-\tau_{\nu}'(1^+)=\liminf-\tau_{\nu,n}'(1)=\liminf\left(-q\tau_{\mu,n}'(q)+\tau_{\mu,n}(q)\right).$$

The following lemma then implies the lower bound.
\begin{Lemma}\label{az} We have $\nu\left(\underline E_{-\tau'(q^-)}\cap\overline F_{-\tau'(q^+)}\right)=1$.
\end{Lemma}
\begin{Remark} Contrary to more regular situations (e.g \cite{bnbh, heurt, ol}), we cannot obtain the more precise result $\nu\left(\overline E_{-\tau'(q^-)}\cap\underline F_{-\tau'(q^+)}\right)=1$ where $$ \overline E_{\alpha}=\left\{x\; ;\;
\alpha(x)\ge \alpha\right\},  \;   \underline F_{\alpha}=\left\{x\; ;\;
\limsup_{n\to\infty}\alpha_n(x)\le \alpha\right\}.$$
\end{Remark}
\begin{Proof}{\bf of Lemma \ref{az}}
For $\eta>0$ we put $\beta=-\tau_{\mu}'(q^-)+\eta$ and we prove that $\nu(\Sigma\setminus\underline E_{\beta})~=~0$. 
In a similar way, it can be shown  that
$\nu(\Sigma\setminus\overline F_{\gamma})=0$ for $\gamma<-\tau_{\mu}'(q^+)$. The lemma then easily follows.

It suffices to show that 
$\displaystyle \Sigma\setminus \underline E_{\beta}=\left\{x\in\Sigma\; ;\; 
\liminf_{n\to\infty}\alpha_{n}(x)> \beta\right\}$ is of 0 $\nu$-measure. 
Consider the collection $\Re_{n}(\beta)$ of cylinders $I\in\Im_{n}$  satisfying
$\ds\frac{\log\mu(I)}{\log|I|}>\beta$. It is clear that $\ds \Sigma\setminus \underline E_{\beta}\subset\liminf_{n\to\infty}\tilde\Re_{n}(\beta)$ with $\tilde\Re_{n}(\beta)=\{x\in \Sigma\;;\;I_n(x)\in \Re_{n}(\beta)\}$.

Let $(\tau_{\mu,n_k})_{k\in\N}$ be the subsequence of $(\tau_{\mu,n})_{n\in\N}$ such that $\lim_{k\to\infty}\tau_{\mu,n_k}(q)=\tau_{\mu}(q)$. 
Using the convergence of $\tau_{\mu,n_k}(q)$ we can choose
(and fix) $t<0$ such that for  $k$ big enough
$$\tau_{\mu}(q+t)-\tau_{\mu,n_k}(q)<-\left(\beta-\frac{\eta}{2}\right)t=\left(\tau_{\mu}'(q^-)-\frac{\eta}{2}\right) t.$$
Since   $\ds \mu(I)^{-t}|I|^{\beta t}\le 1$ if  $I\in \Re_{n}(\beta)$, we have 
\begin{eqnarray*}
\nu(\tilde\Re_{n_k}(\beta))=\sum_{I\in\Re_{n_k}(\beta)}\nu(I)&=&\sum_{I\in\Re_{n_k}(\beta)}\mu(I)^q|I|^{\tau_{\mu,n_k}(q)}=
\sum_{I\in\Re_{n_k}(\beta)}\mu(I)^{q+t}|I|^{\tau_{\mu,n_k}(q)-\beta t}\mu(I)^{-t}|I|^{\beta t}\\
&\le&\sum_{I\in\Re_{n_k}(\beta)}\mu(I)^{q+t}|I|^{\tau_{\mu,n_k}(q)-\beta t}\le 
|I|^{-\frac{\eta}{4}t}\sum_{I\in\Im_{n_k}}\mu(I)^{q+t}|I|^{\tau_{\mu}(q+t)-\frac{\eta}{4}t}\\
&\le& |I|^{-\frac{\eta}{4}t} \sum_{I\in\Im_{n_k}}\mu(I)^{q+t}|I|^{\tau_{\mu,n_k}(q+t)}=|I|^{-\frac{\eta}{4}t}.
\end{eqnarray*}
For the last inequality,  we used the fact that $\ds\tau_{\mu}(q+t)=\limsup\tau_{\mu,n}(q+t)$. We deduce that  $$\liminf_{n\to\infty}\nu(\tilde\Re_{n}(\beta))=0$$ 
 and the lemma easily follows.

%Fix $\delta>0$ and  $n_0\in\N$ such that for all $n\ge n_0$, $\tau_{\mu}(q)<\tau_n(q)+\delta/2$. We consider, 
%for every $x\in \Sigma$ the bigest cylinder $\displaystyle J(x)\in\cup_{n\ge n_0}\Im_n$
%such that $\frac{\log\mu(J(x))}{\log|J(x)|}\ge\alpha+\delta$. The collection of these cylinders forms
%a covering $\Re$ of $\underline E_{\alpha}$ and we have
\end{Proof}
%$$\sum_{J\in\Re}\nu(J)\le\frac1C \sum_{J\in\Re}\mu(J)^q|J|^{\tau_{\mu}(q)-\delta/2}\le \frac1C
%\sum_{J\in\Re}|J|^{q(\alpha-\delta)}|J|^{\tau_{\mu}(q)+\delta/2}=
%\frac1C\sum_{J\in\Re}|J|^{q\alpha+\tau_{\mu}(q)+\delta/2-q\delta}.$$
The proof of Theorem \ref{1} is now completed.\end{Proof}

%Let $f$ and $g$ be the functions defined by $f(t)=\dim \underline E_t$ and $g(t)=\dim \overline F_t$ . 
%Obviously, $f$ is increasing and $g$ is decreasing. 
%Recall that $t$ is a non-stationary point of a monotone function $h$ if $h(s)\not=h(t)$ for all $s\not=t$.

%Since $E_{\alpha}=\underline E_{\alpha}\setminus\bigcup_{\beta< \alpha} \underline E_{\beta}$,  we deduce from theorem \ref{1} the following.
%\begin{Remark}
%If $\alpha=-\tau'(q^-)$ for $q>0$  is a non-stationary point of $f$ there a sequence of $q_m\le q$ such that $\alpha_m=-\tau'(q_m^-)$  
%are non-stationary points of $f$ converging to $\alpha$ and
%$$\liminf_{n\to\infty}-q_m\tau_{\mu,n}'(q_m)+\tau_{\mu,n}(q_m)\le \dim E_{\alpha_m}=\dim \underline E_{\alpha_m} \le \tau^*(\alpha_m).$$
%If $\alpha=-\tau'(q^+)$ for $q<0$  is a non-stationary point of $g$ there a sequence of $q_m\ge q$ such that $\alpha_m=-\tau'(q_m^+)$  
%are non-stationary points of $g$ converging to $\alpha$ and
%$$\liminf_{n\to\infty}-q_m\tau_{\mu,n}'(q_m)+\tau_{\mu,n}(q_m)\le \dim F_{\alpha_m}=\dim \overline F_{\alpha_m} \le \tau^*(\alpha_m).$$

%We conjecture  that under the same conditions on $\alpha$ we should also have $\dim E_{\alpha}=\dim \underline E_{\alpha}$ 
%($\dim F_{\alpha}=\dim \overline F_{\alpha}$ respectively).
%
%\end{Remark}

\section{Proof of Theorem \ref{4}}
%Some conditions  ensuring  the validity of  multifractal formalism}
 We will use  the following result.

\begin{Propo}\label{2}
For  $q\in\R$, let   $(\tau_{\mu,n_k})$ be 
the subsequence of $(\tau_{\mu,n})$ such that 
$$\displaystyle \lim_{k\to\infty}\tau_{\mu,n_k}(q)=\limsup_{n\to\infty}\tau_{\mu,n}(q)=\tau_\mu(q).$$ Then, we have $$\tau_{\mu}'(q^-)\le\liminf_{k\to\infty}\tau_{\mu,n_k}'(q)\le \limsup_{k\to\infty}\tau_{\mu,n_k}'(q)\le \tau_{\mu}'(q^+)$$ where $\tau_{\mu}'(q^-)$ and $\tau_{\mu}'(q^+)$ stand  for the left and the right hand d\'erivative   of $\tau_{\mu}$ at $q$.

 Hence, if $\tau_{\mu}'(q)$ exists,  we have  
$$\lim_{k\to\infty}\tau_{\mu,n_k}'=\tau_{\mu}'(q).$$
\end{Propo}

\begin{Proof} 
 We only prove the inequality $\displaystyle\limsup_{k\to\infty}\tau_{\mu,n_k}'(q)\le \tau_{\mu}'(q^+)$. The proof of   $\displaystyle\tau_{\mu}'(q^-)\le\liminf_{k\to\infty}\tau_{\mu,n_k}'(q)$   is similar. 

 Take   $\epsilon>0$ and   $\tilde q>q$ such that 
$$\left|\frac{\tau_{\mu}(\tilde q)-\tau_{\mu}(q)}{\tilde q-q}-\tau_{\mu}'(q^+)\right|<\epsilon/3.$$
 %We  consider $(\tilde n_k)$ such that 
 %$\displaystyle \lim_{k\to\infty}\tau_{\mu,\tilde n_k}(\tilde q)=
%\tau_{\mu}(\tilde q)$.

We can chose $k$ big enough to have 
\begin{eqnarray*}\label{restrictions2}
\frac{|\tau_{\mu,n_k}(q)-\tau_{\mu}(q)|}{|\tilde q-q|}&<&\epsilon/3\end{eqnarray*} and 
%\frac{|\tau_{\mu,\tilde n_k}(\tilde q)-\tau_{\mu}(\tilde q)|}{|\tilde q-q|}&<&\epsilon/8\\
\begin{eqnarray*} \tau_{\mu,n_k}(\tilde q)&\le& \tau_{\mu}(\tilde q)+(\tilde q-q)\epsilon/3.
\end{eqnarray*}
We then obtain
\begin{eqnarray*}
\tau_{\mu}'(q^+) &\ge& \displaystyle \frac{\tau_{\mu}(\tilde q)-\tau_{\mu}(q)}{\tilde q-q}-\epsilon/3\\
&=& \frac{\tau_{\mu}(\tilde q)-\tau_{\mu, n_k}(\tilde q)+\tau_{\mu,n_k}(\tilde q)-\tau_{\mu,n_k}(q)+\tau_{\mu,n_k}(q)-\tau_{\mu}(q)}{\tilde q-q}-\epsilon/3\\
&\ge& -\epsilon/3+  \tau_{\mu,n_k}'(q) -\epsilon/3-\epsilon/3= \tau_{\mu,n_k}'(q)- \epsilon
\end{eqnarray*}
and the proof easily follows.
\end{Proof}
%In particular we have showed the following
%\begin{Corollary}
%Let $q$ be a differentiability point for $\tau_{\mu}$. For all $\epsilon>0$ there exists $\delta>0$ such that
%$|\tilde q-q|<\delta$ implies $|\tau_{\mu,\tilde n_k}(\tilde q)-\tau_{\mu,n_k}(q)|<\epsilon$, where
%$(\tilde n_k)$ and $(n_k)$ are such that 
 %$\displaystyle \lim_{k\to\infty}\tau_{\mu,\tilde n_k}(\tilde q)=
%\limsup_{n\to\infty}\tau_{\mu,n}(\tilde q)$,
 %$\displaystyle \lim_{k\to\infty}\tau_{\mu, n_k}(q)=
%\limsup_{n\to\infty}\tau_{\mu,n}(q)$ and $k=k(q,\tilde q,\epsilon)$ big enough.
%\end{Corollary}

We can now  prove Theorem \ref{4}.
 
\begin{Proof}{\bf of Theorem \ref{4}}.
Let $\nu$ be the Gibbs-measure defined in Lemma \ref{nuthing}. Since 
$$\tau_{\nu,n}(s)=\tau_{\mu,n}(qs)-s\tau_{\mu,n}(q)$$ we get 
$$\tau_{\nu,n}'(1)=q\tau_{\mu,n}'(q)-\tau_{\mu,n}(q).$$
Using the convergence of $\tau_{\mu,n}(q)$ we deduce from Proposition \ref{2} that
$$\lim_{n\to\infty}\tau_{\nu,n}'(1)=\lim_{n\to\infty}\left(q\tau_{\mu,n}'(q)-\tau_{\mu,n}(q)\right)=q\tau_\mu'(q)-\tau_\mu(q).$$
Since $\nu$ is also an inhomogeneous Bernoulli product, we deduce from  remark \ref{Heurt} that $\tau_{\nu}'(1)$ exists and 
$$\dim\nu=\Dim\nu=-\tau_{\nu}'(1)=-q\tau_\mu'(q)+\tau_\mu(q).$$

On the other hand, for $I\in\Im_n$, we have 
$$\frac{\log\nu(I)}{\log |I|}=q\frac{\log\mu(I)}{\log|I|}+\tau_{\mu,n}(q).$$
 Since 
 $$\lim_{n\to\infty}\frac{\log\nu(I_n(x))}{\log |I_n(x)|}=\dim\nu=\Dim\nu\; \,;\nu\mbox{-a.s.}$$
 we obtain that $\ds \lim_{n\to\infty}\frac{\log\mu(I_n(x))}{\log |I_n(x)|}=-\tau_\mu'(q)$, $\nu$-a.s.
 We conclude that $$\dim\left(E_{\alpha}\cap F_{\alpha}\right)\ge \dim\nu=\tau_\mu^*(\alpha).$$
 The opposite inequality being always valid, the proof  is done.
\end{Proof}

We end this section with a few comments about Theorem \ref{4}. 

  As mentionned in the introduction of the paper,  the validity of the multifractal formalism is often easier  to obtain for the values of $\alpha$ that can be written $\alpha=-\tau'(q)$. However, the following example shows that there exist inhomogeneous Bernoulli products that do not satisfy the multifractal formalism even at their differentiability points  $\alpha=-\tau'(q)$. Thus, the convergence of the sequence  $(\tau_{n}(q))$  is really necessary for the validity of the multifractal formalism in our context.

To see that, let $(n_k)_{k\ge 1}$ be a sequence of integers such that $$n_1=1,\quad n_k<n_{k+1} \quad \text{and} \quad \lim_{k\to +\infty} \frac{n_{k+1}}{n_k}=+\infty,$$  and  consider the inhomogeneous Bernoulli product $\mu$ given by the sequence $(p_n)$ such that $$p_i=p \quad \text{if} \quad n_{2n-1}\le i < n_{2n}\quad \text{and} \quad p_i=\tilde p \quad \text{if} \quad n_{2n}\le i < n_{2n+1},$$ with $0<p<\tilde p< 1/2$.

The calculation of the function $\tau$ is classical. By observing that $$\mu(I_{\epsilon_1...\epsilon_n0})^q +\mu(I_{\epsilon_1...\epsilon_n1})^q=[(p_{n+1}^q+(1-p_{n+1})^q]\mu(I_{\epsilon_1...\epsilon_n})^q,$$  we easily deduce that $$\sum_{I\in \mathcal{F}_n}\mu(I)^q=\prod_{k=1}^n[p_k^q+(1-p_k)^q].$$ Then, if $N_n$ is the number of integer $k\le n$ such that $p_k=p$, we have $$\tau_n(q)=\frac{N_n}{n} \log_2(p^q+(1-p)^q)+(1-\frac{N_n}{n})\log_2(\tilde p^q+(1-\tilde p)^q).$$ Since that $\liminf_n \frac{N_n}{n}=0$ and  $\limsup_n \frac{N_n}{n}=1$, we get  $$\tau(q)=\sup(\log_2(p^q+(1-p)^q),  \log_2(\tilde p^q+(1-\tilde p)^q).$$  So, except for $q=0$ and $q=1$, $\tau'(q)$ exists. Moreover, $$\forall I\in \mathcal{F}_n,\quad \mu(I)\ge \mu(I_{00\cdots0})=p^{N_n}\tilde p^{n-N_n}.$$ Thus, $$ \forall I\in \mathcal{F}_n,\quad -\frac{\log_2(\mu(I))}{n}\le \frac{N_n}{n} (-\log_2p)+(1-\frac{N_n}{n})(-\log_2 \tilde p),$$ and we have  $$\forall I\in \mathcal{F}_n,\quad \liminf_n  -\frac{\log_2(\mu(I))}{n}\le \inf(-\log_2  p,-\log_2 \tilde p)=-\log_2 \tilde p.$$ Finally,  if $-\log_2 \tilde p<\alpha= -\tau'(q)<-\log_2  p$, we have $E_{-\tau'(q)}=\emptyset$ and the multifractal formalism is not satisfied for a such $\alpha$.   

Moreover, this example shows that the function $\tau$ may be not differentiable at the positive values of $q$. Therefore, the situation differs from this one obtained by Heurteaux \cite{heurt} for quasi-Bernoulli measure for which $\tau$ is differentiable on $\R$. It also differs from this one  obtained by Testud in \cite{testud}   for weak quasi-Bernoulli  measure for which    the phase transitions only appears for $q<0$. 

In fact, the function $\tau$ of an  inhomogeneous Bernoulli product may be very irregular. This is the object following section.

\section{Proof of  Theorem \ref{3}}
 From now, we denote par $\tau(p,.)$ the $\tau$ function of the  homogeneous Bernoulli product of parameter $p$. Moreover, whenever we use the notation $p_i$ for a weight in $(0,1)$ we will also note $\tau_i=\tau(p_i,.)$.

Before the proof of  Theorem \ref{3}, we present a few lemmas.

\begin{Lemma}\label{subsidiary} For any $p_1<p_2<p_3$  in  $(0,1/2)$   consider  the functions $\tau_1=\tau(p_1,.),\tau_2=\tau(p_2,.)$ and $\tau_3=\tau(p_3,.)$. We have that 
$\ds \frac{\tau_1-\tau_2}{\tau_2-\tau_3}$ is decreasing on $(1,+\infty)$.
\end{Lemma}
Although the proof only uses elementary calculus, it is a litte bit ``tricky'' and cannot be omitted.
\begin{Proof}{\bf of Lemma \ref{subsidiary}} Taking into account the trivial equality
$$\tau(p',q)-\tau(p'',q)=\int_{p''}^{p'}\frac{\partial\tau}{\partial p}(p,q)dp$$ we only need to show that if $p'<p''$ then\; \;
$\ds\frac{\partial\tau}{\partial p}(p',q):\frac{\partial\tau}{\partial p}(p'',q)$\; \;is decreasing on $q\in(1,\infty)$.
We get
\begin{eqnarray*}
\ds\frac{\partial\tau}{\partial p}(p',q):\frac{\partial\tau}{\partial p}(p'',q) &=& \frac{1}{p'}\frac{1-(-1+1/p')^{q-1}}{1+(-1+1/p')^{q}}:\frac{1}{p''}\frac{1-(-1+1/p'')^{q-1}}{1+(-1+1/p'')^{q}}\\
&=& p''\frac{1-{s_1}^{q-1}}{1+{s_1}^{q}}:p'\frac{1-{s_2}^{q-1}}{1+{s_2}^{q}}
\end{eqnarray*}
where $s_1=-1+1/p'>1$ and $s_2=-1+1/p''>1$.

If we set $\ds f(s,q)=\ln\frac{1-{s}^{q-1}}{1+{s}^{q}}$, with $s,q>1$, it is  sufficient to prove that $\ds\frac{\partial f}{\partial s}f(s,q)$ is decreasing in $q$.
We calculate
$$\ds\frac{\partial f}{\partial s}f(s,q)=\frac{(q-1)s^{q-2}}{s^{q-1}-1}-\frac{qs^{q-1}}{s^{q}+1}.$$
By  multiplying  by $s$, we  need to show that
$\ds \frac{(q-1)s^{q-1}}{s^{q-1}-1}-\frac{qs^{q}}{s^{q}+1}$ is decreasing  which is equivalent to 
$q-1+\frac{q-1}{s^{q-1}-1}-q+\frac{q}{s^{q}+1}$ being decreasing. 

Put $Q=q-1$ ; it remains to show that
$\ds \frac{q-1}{s^{q-1}-1}+\frac{q}{s^{q}+1}=\frac{Q}{s^Q-1}+\frac{Q}{s^{Q+1}+1}+\frac{1}{s^{Q+1}+1}$ decreases in $Q>0$.  The last term being decreasing it suffices to show that $\ds\frac{Q}{s^Q-1}+\frac{Q}{s^{Q+1}+1}$ is doing the same.
By taking derivatives we need to show that 
$$s^{Q+1}(s^Q-1-s^Q\ln s^Q)+s^Q-1-\ln s^Q\le 0.$$ Since, $(s^Q-1-s^Q\ln s^Q)<0$, it suffices to show that $$ s^{Q}(s^Q-1-s^Q\ln s^Q)+s^Q-1-\ln s^Q=s^{2Q}-s^{2Q}\ln s^Q-\ln s^Q-1=g(s^Q)\le 0$$ where $g(x)=x^2-x^2\ln x-\ln x-1.$ Moreover,   the sign of  $g'(x)=x-x\ln x^2-1/x$ is  the same   of the sign of     $x^2-x^2\ln x^2-1$ if $x>1$. Since, $y-1 \le y\ln y$ for $y>1$, we deduce that $g$ is decreasing on $(1,+\infty)$. By observing than  $g(1)=0$, we obtain that $\ds\frac{\partial f}{\partial s}f(s,q)$ is decreasing  on $(1,+\infty)$ and the    Lemma \ref{subsidiary} is proved.
\end{Proof}

\begin{Lemma}\label{basic}
Take $\tau=\lambda\tau(p_1,.)+(1-\lambda)\tau(p_2,.)$ with $0<p_1<p_2<1/2$ and $\lambda\in(0,1)$. For $p_0\in(0,1/2)$
one of the following occurs:
\begin{enumerate}
\item  either $\tau(q)\not= \tau(p_0,q)$, for all $q>1$,
\item either there exists $q_0>1$ such that $\tau(q)>\tau(p_0,q)$ for $1<q<q_0$ and $\tau(q)< \tau(p_,q)$ for $q>q_0$. In this case  $q_0$ is then the unique point of $(1,+\infty)$  for which $\tau(q)= \tau(p_0,q)$.
\end{enumerate}
\end{Lemma}

\begin{Proof}{\bf of Lemma \ref{basic}.}
Let us first remark that  $\tau$ and $\tau(p_0,.)$ can coincide at one point  only if $p_0\in(p_1,p_2)$. Moreover,
$\tau(q)=\tau(p_0,q)$ implies
$$\frac{\tau(p_1,q)-\tau(p_0,q)}{\tau(p_0,q)-\tau(p_2,q)}=\frac{1-\lambda}{\lambda}.$$
By Lemma \ref{subsidiary} this can only occur once and Lemma \ref{basic} easily follows on the decreasing property of the ratio.
\end{Proof}

\begin{Lemma}\label{system} 
Take $\lambda_1,\lambda_2\in(0,1)$ such that $\lambda_1+\lambda_2=1$ , $1<p_1<p_2<1/2$ and 
set $\tau=\lambda_1\tau_1+\lambda_2\tau_2$. Fix $1<q_1<q_2<+\infty$ and consider $p_1<p_4<p_2<p_5<1/2$ such that $\tau(p_4,q_1)=\tau(q_1)$. Then there is a unique convex combination $\tilde \tau$ of $\tau_1,\tau_4$ and $\tau_5$ such that
$$\tilde\tau(q_1)=\tau(q_1)\mbox{ and }\tilde\tau(q_2)=\tau(q_2).$$
Furthermore, for $i=1,2$, we have $\tau'(q_i)\not=\tilde\tau'(q_i)$ and $\tau(q)\not=\tilde\tau(q)$ if $1<q\not=q_i$.
\end{Lemma}
\begin{Proof}{\bf of Lemma \ref{system}.}

 First note that it is easy to see that there exists $p_4\in (p_1,p_2)$ such that $\tau(p_4,q_1)=\tau(q_1)$.

It then suffices to  show that the linear system 
\begin{equation*}\left\{\begin{array}{llllcl}
&\lambda_3\tau_1(q_1)&+\lambda_4\tau_4(q_1)&+\lambda_5\tau_5(q_1)&=&\tau(q_1)\\
&\lambda_3\tau_1(q_2)&+\lambda_4\tau_4(q_2)&+\lambda_5\tau_5(q_2)&=&\tau(q_2)\\
&\lambda_3&+\lambda_4&+\lambda_5&=&1
\end{array} \right.\hspace{1cm}{ \bf (S)}\end{equation*}
has a unique positive solution $(\lambda_3,\lambda_4,\lambda_5)$.  The existence of a unique solution is easy to verify. Let us show that this solution is positive.

First note that  $\lambda_4\not = 1$. Indeed, if   $\lambda_4= 1$, since $\tau(q_1)=\tau_4(q_1)$, we have $\lambda_3(\tau_1(q_1)-\tau_5(q_1))=0$. Thus, $ \lambda_3=\lambda_5=0$ and $\tau(q_2)=\tau_4(q_2)$ which is not possible by Lemma \ref{basic}. 

Therefore, since $\tau(q_1)=\tau_4(q_1)$, the first equation of the system gives that 

\begin{eqnarray}\label{s}\frac{\lambda_3}{\lambda_3+\lambda_5}\tau_1(q_1)+\frac{\lambda_5}{\lambda_3+\lambda_5}\tau_5(q_1)=\lambda_1\tau_1(q_1)+\lambda_2\tau_2(q_1)\in (\tau_5(q_1), \tau_1(q_1))\end{eqnarray}
This implies   that $\frac{\lambda_3}{\lambda_3+\lambda_5}\in (0,1)$. We deduce that   $\lambda_3 \lambda_5>0$.  Moreover,  since $\tau_5<\tau_2$, we also  have $\frac{\lambda_3}{\lambda_3+\lambda_5} >\lambda_1$. 

Let us show that $\lambda_3$ and $\lambda_5$ are positive. Otherwise,  by the above remark, we have $\lambda_3<0$, $\lambda_5<0$ and $\lambda_4> 0$. By the system (S) we have 
\begin{eqnarray*}\tau_4(q)=\frac{\lambda_1-\lambda_3}{\lambda_4}\tau_1(q)+\frac{\lambda_2}{\lambda_4}\tau_2(q)-\frac{\lambda_5}{\lambda_4}\tau_5(q)\end{eqnarray*} at the points $q=q_1$ and $q=q_2$. We then obtain that  $$\frac{\lambda_1-\lambda_3}{\lambda_4}\frac{\tau_1-\tau_4}{\tau_4-\tau_2}(q)=\frac{\lambda_2}{\lambda_4}-\frac{\lambda_5}{\lambda_4}\frac{\tau_4-\tau_5}{\tau_4-\tau_2}(q)$$ for $q=q_1$ and $q=q_2$. Since $p_1<p_4<p_2$, by Lemma \ref{subsidiary} the function $\frac{\tau_1-\tau_4}{\tau_4-\tau_2}$ is decreasing. On the other hand,  since  $p_4<p_2<p_5$,  Lemma \ref{subsidiary}  implies that the function $\frac{\tau_4-\tau_5}{\tau_4-\tau_2}=1+\frac{\tau_2-\tau_5}{\tau_4-\tau_2}$ is increasing. Thus, these functions cannot coincide at two  points so   we conclude  that $\lambda_3$ and $\lambda_5$ are positive.

Let us now prove that $\lambda_4>0$. By \eqref{s} we have $$\frac{\lambda_3}{\lambda_3+\lambda_5}\tau_1(q_1)+\frac{\lambda_5}{\lambda_3+\lambda_5}\tau_5(q_1)=\lambda_1\tau_1(q_1)+\lambda_2\tau_2(q_1)$$ which gives that $$\lambda_2\tau_2(q_1)=\left(\frac{\lambda_3}{\lambda_3+\lambda_5}-\lambda_1\right)\tau_1(q_1)+\frac{\lambda_5}{\lambda_3+\lambda_5}\tau_5(q_1).$$ Using Lemma \ref{subsidiary}, for $q> q_1$   we get
$$\lambda_2\tau_2(q)> \left(\frac{\lambda_3}{\lambda_3+\lambda_5}-\lambda_1\right)\tau_1(q)+\frac{\lambda_5}{\lambda_3+\lambda_5}\tau_5(q)$$ and $$\lambda_1\tau_1(q)+\lambda_2\tau_2(q)> \frac{\lambda_3}{\lambda_3+\lambda_5}\tau_1(q)+\frac{\lambda_5}{\lambda_3+\lambda_5}\tau_5(q).$$ In particular, for $q=q_2$ we find that $$\lambda_3\tau_1(q_2)+\lambda_5\tau_5(q_2)+\lambda_4\tau(q_2)<\tau(q_2)=\lambda_3\tau_1(q_2)+\lambda_4\tau_4(q_2)+\lambda_5\tau_5(q_2)$$ and we deduce that $$\lambda_4 \tau(q_2)<\lambda_4 \tau_4(q_2).$$ Since $\tau(q_1)=\tau_4(q_1)$, it follows from Lemma \ref{subsidiary} that $\lambda_4>0$.

The last assertion follows directly from the independancy of the vector families
$$\left\{\left(\begin{array}{c}\tau_1(q_1)\\ \tau_4(q_1)\\ \tau_5(q_1)\end{array}\right), \left(\begin{array}{c}\tau_1(q_2)\\ \tau_4(q_2)\\ \tau_5(q_2)\end{array}\right), \left(\begin{array}{c}\tau_1'(q_i)\\ \tau_4'(q_i)\\ \tau_5'(q_i)\end{array}\right)\right\} $$ and
$$\left\{\left(\begin{array}{c}\tau_1(q_1)\\ \tau_4(q_1)\\ \tau_5(q_1)\end{array}\right), \left(\begin{array}{c}\tau_1(q_2)\\ \tau_4(q_2)\\ \tau_5(q_2)\end{array}\right), \left(\begin{array}{c}\tau_1(q)\\ \tau_4(q)\\ \tau_5(q)\end{array}\right)\right\} ,$$
which can be easily established.
\end{Proof}

\begin{Remark}\label{proximity}
In the proof of Lemma \ref{system} it is clear that when $p_5$ is close to $p_2$, the solution of the system ({\bf S}) converges to $(\lambda_1,0,\lambda_2)$  and $\tilde\tau$ converges to $\tau$.
\end{Remark}

The following result generalizes Lemma \ref{system} for any convex combination of functions $\tau(p_i,.)$.
 
\begin{Lemma}\label{cocorico} Let $\tau$ be a convex combination of functions $\tau(p_i,.)$ where $0<p_i\le 1/2$, $i=1,...,n\ge 2$. For any $1<q_1<q_2<\infty$ there exists another convex combination $\tilde \tau$ of functions $\tau(p_j',.)$ such that \begin{itemize}
\item $\tilde\tau(q_i)=\tau(q_i)$ and $\tilde\tau'(q_i)\not=\tau'(q_i)$, $i=1,2$,
\item for $q\not\in\{q_1,q_2\}$, $\tilde\tau(q)\not=\tau(q)$.
\end{itemize}
\end{Lemma}
\begin{Proof}{\bf of Lemma \ref{cocorico}.}

 The  case $n=2$ is given by Lemma \ref{system}. The case $n>2$ is easy to derive. Suppose   $\tau=\sum_{k=1}^{n}\lambda_k\tau(p_k,.)$ and let $\tau_1=\tau(p_1,.)$ and $\tau_2=\tau(p_2,.)$ be the first two functions of the convex combination.   By  Lemma \ref{system} there exists a convex combination $\hat\tau$ of three $\tau(p,.)$ functions such that 
\begin{enumerate}
\item $\frac{1}{\lambda_1+\lambda_2}\left(\lambda_1\tau_1(q_i)+\lambda_2\tau_2(q_i)\right)=\hat\tau(q_i)$ , for $i=1,2$  
\item $\frac{1}{\lambda_1+\lambda_2}\left(\lambda_1\tau_1'(q_i)+\lambda_2\tau_2'(q_i)\right)\not =\hat\tau'(q_i)$ for $i=1,2$.   
\end{enumerate}
The function $\tilde\tau=(\lambda_1+\lambda_2)\hat\tau+\sum_{k=3}^{n}\lambda_k\tau(p_k,.)$ satisfies then the conclusion of Lemma~\ref{cocorico}.\end{Proof}

The following lemma is easy and the proof is  left to the reader.

\begin{Lemma}\label{previous}
For any $p_1,...,p_n$ and any convex combination $\tau$ of $\tau(p_1,.),...,\tau(p_n,.)$ there exists an inhomogeneous Bernoulli measure $\mu$ whose multifractal spectrum equals $\tau$.
\end{Lemma}

We can now prove Theorem \ref{3}.

\begin{Proof}{\bf of Theorem \ref{3}} 

 In fact, and in order to avoid technicalities, we only prove the following easier version of Theorem~\ref{3} and then indicate the changes needed to extend the proof in the general case.

\begin{Theorem}\label{pro}
 There exists an inhomogeneous Bernoulli product $\mu$ such that the spectrum $\tau$ of $\mu$ has an infinite set of  the phase transitions on $(1, +\infty)$. Moreover, this set has a finite point of accumulation. 
\end{Theorem}
\begin{Proof}
The strategy of the demonstration  is the following : we first find inhomogeneous Bernoulli products that are not derivable at a finite number of predefined points and we construct the measure $\mu$ using Cantor's diagonal argument.

Fix $(q_n)_{n\ge 1}$ a sequence of real numbers  nested in the sense that $1<q_1<...<q_{2n-1}<q_{2n+1}<q_{2n+2}<q_{2n}<...<q_2$  for all $n\ge 1$ and $\cap_n(q_{2n+1}, q_{2n+2})=\{q_0\}$. In particular, $\lim q_n=q_0$.  Let $p_1,p_2\in (0,1/2)$ and consider $\tau_1=\frac12\tau(p_1,.)+\frac12\tau(p_2,.)$. By Lemma~\ref{previous} we can construct a Bernoulli product $\mu_1$ of spectrum $\tau_1$. Then, Lemma~\ref{cocorico} implies  the existence of a convex combination $\tau_2$ of $\tau(p_i,.)$'s functions, such that 

$$\tau_1(q_i)=\tau_2(q_i)\mbox{ , for }i=1,2\mbox{  and }\tau_1'(q_i)\not=\tau_2'(q_i).$$

We can define a measure $\mu_2$ of spectrum $\tau_2$. Using $\mu_1$ and $\mu_2$, we can construct a measure $\nu_2$ of spectrum $\rho_2=\max\{\tau_1,\tau_2\}$. To do that,  we take  a sequence of integers $(\ell_k)_k$ such that $\displaystyle\frac{\ell_{k+1}}{\sum_1^k\ell_i}\to\infty$. On dyadique intervals of length between $2^{-\ell_{2k}}$ and $2^{-\ell_{2k+1}}$ we apply the weight distribution of $\mu_1$ and on dyadique intervals of length between $2^{-\ell_{2k+1}}$ and $2^{-\ell_{2k+2}}$ we apply the weight distribution of $\mu_2$, where $k\in \N$. It is easy to verify that the resulting inhomogeneous measure $\nu_2$ has spectrum $\rho_2=\max\{\tau_1,\tau_2\}$. The spectrum of $\nu_2$ is not differentiable at $q_1$ and $q_2$. 

We proceed by induction to construct a measure $\nu_n$ which has a non differentiable spectrum for points $q_1,\cdots,q_{2n-2}$.  Suppose the measures $\nu_1=\mu_1$, $\mu_2$, $\nu_2$,..., $\mu_n$, $\nu_n$ constructed and denote by $\rho_n= \max\{\tau_1,...\tau_n\}$ where $\tau_i$ is  the spectrum of the measure $\mu_i$, $i\in\{1...,n\}$. %We assume that that on every interval $[q_{2i+1},q_{2i+2}]$ , where $i\le n$, the spectrum of $\nu_n$ equals  $\max\{\tau_1,...\tau_n\}$ and is realized by $\tau_i$. 
Let us construct $\mu_{n+1}$ and $\nu_{n+1}$. 

 One of the following two cases hold:

\begin{itemize}
\item[{\bf Case 1}]
Lemma~\ref{cocorico} provides a function $\tau_{n+1}$ satisfying : 
\begin{enumerate}
\item $\tau_{n+1}(q_{2n-i})=\rho_n(q_{2n-i})$  for $i=0,1$ and $\tau_{n+1}(q)\not=\rho_n(q)$ if $q\not \in \{q_{2n-1},q_{2n}\}$,
\item $\tau_{n+1}'(q_{2n-1})>\rho_n'(q_{2n-1})$ \; , \; $\tau_{n+1}'(q_{2n})<\rho_n'(q_{2n})$  
\end{enumerate}
Therefore  we have $\tau_{n+1}>\rho_n \mbox{ on }(q_{2n-1},q_{2n})$ and 
$\tau_{n+1}<\rho_n$ on $(1,\infty)\setminus [q_{2n-1},q_{2n}].$

Let $\mu_{n+1}$ be the inhomogeneous Bernoulli measure of spectrum $\tau_{n+1}$.  To define the measure $\nu_{n+1}$ we use the previous procedure convenably adapted:
Take $(\ell_k)_k$ a sequence of integers such that $\displaystyle\frac{\ell_{k+1}}{\sum_1^k\ell_i}\to\infty$. On dyadique intervals of length between $2^{-\ell_{(n+1)k+i}}$ and $2^{-\ell_{(n+1)k+i+1}}$ apply the weight distribution of $\mu_i$, where $i=1,...,n+1$ and $k\in \N$. It is easy to verify that the resulting inhomogeneous measure $\nu_{n+1}$ has spectrum $\rho_{n+1} =\max\{\tau_1,...\tau_{n+1}\}$ on $(1,\infty)$. Remark that this spectrum equals $\tau_{n+1}$ on $[q_{2n-1},q_{2n}]$ and $\rho_n=\max\{\tau_1,...\tau_n\}$ elsewhere on $[1,\infty)$. Clearly, in this case, the function $\rho_{n+1}= \max(\tau, \tau_{n+1})$ is not differentiable at $q_1,\cdots,q_{2n}$.

\

\item[{\bf Case 2}] Lemma~\ref{cocorico} provides for all choices of $p_5>p_2$ a function $\tau_{n+1}$ satisfying: 
\begin{enumerate}
\item $\tau_{n+1}(q_{2n-i})=\rho_n(q_{2n-i})$ , for $i=0,1$ and $\tau_{n+1}(q)\not=\rho_n(q)$ if $q\not \in \{q_{2n-1},q_{2n}\}$,
\item $\tau_{n+1}'(q_{2n-1})<\rho_n'(q_{2n-1})$ \; , \; $\tau_{n+1}'(q_{2n})>\rho_n'(q_{2n})$  
\end{enumerate}
In this case,
$$\tau_{n+1}<\rho_n \mbox{ on }(q_{2n-1},q_{2n})\mbox{ and }
\tau_{n+1}>\rho_n\mbox{ on }(q_{2n-3},q_{2n-1})\cup(q_{2n},q_{2n-2}).$$
The function $\tilde \rho_{n+1}= \max(\rho_n, \tau_{n+1})$ is not differentiable at $q_{2n-1}$ and $q_{2n}$ but we lose the phase transitions $q_{2n-3}, q_{2n-2}$ and we don't know what happens for the  other phase transitions  $q_1,\cdots,q_{2n-4}$. 

To avoid this problem we use  remark \ref{proximity}. From this, when $p_5$ converges to $p_2$, $\tau_{n+1}$ converges to the    convex combination $T$ of $\tau(p_i,.)$ functions  which is equal to  $\rho_n$ on $(q_{2n-3},q_{2n-2})$.  Since  $T$ differs from $\rho_n$ on $(q_{2n-5},q_{2n-3})$ and $(q_{2n-2},q_{2n-4})$,  we can choose $p_5$ sufficiently close to $p_2$ such that 
$$\tau_{n+1}\left(\frac{q_{2n-3}+q_{2n-5}}{2}\right)<\rho_n\left(\frac{q_{2n-3}+q_{2n-5}}{2}\right) $$ and $$\tau_{n+1}\left(\frac{q_{2n-2}+q_{2n-4}}{2}\right)<\rho_n\left(\frac{q_{2n-2}+q_{2n-4}}{2}\right) .$$
 We deduce that there exist $q_{2n-5}<q'<q_{2n-3}$ and $q_{2n-2}<q''<q_{2n-4}$ such that $\tau_{n+1}=\rho_n $ at $q'$ and $q''$ and $\tau_{n+1}<\rho_n$ on $(q_1,q')$ and $(q'',q_2)$.

The modified family of $\tilde q_i$'s defined by $$ \tilde q_i =
 \begin{cases}
 q_i & \mbox{if}\;\; i\not\in  \{2n-3,2n-2\}\\ q' &  \mbox{if} \;\; i=2n-3\\  q'' &  \mbox{if} \;\; i=2n-2
 \end{cases}$$  have the same properties as the initial $q_i$'s. Moreover,  $\rho_{n+1}= \max(\rho_n, \tau_{n+1})$ is not differentiable at points $q_1,\cdots,q_{2n}$. We proceed as above for the construction of  the measures $\mu_{n+1}$ and $\nu_{n+1}$ which have  spectra  $\tau_{n+1}$ and $\rho_{n+1}$ respectively.
\end{itemize}

To end the proof we use Cantor's diagonal argument: take $(\ell_k)_k$ as before and  define the measure $\nu$ by applying the weight distribution of $\nu_k$ on dyadique intervals of length between $2^{-\ell_{k}}$ and $2^{-\ell_{k+1}}$. The spectrum of the measure $\nu$ equals then $\tau=\sup_{n\in\N}{\rho_n}=\sup_{n\in\N}{\tau_n}$. By  construction, the set of non-derivability points of the function $\tau$ is infinite and has $q_0$ as accumulation point. \end{Proof} 

\begin{Remark}
The second case of the proof  of Theorem~\ref{pro} seems to be inexistent (in our numerical simulations) but we have not been able to prove that only the first  case  arises.\end{Remark}

Let us now give some hints concerning the proof of Theorem~\ref{3}.

 Fix $(q_n)_n$ a sequence of real numbers, dense in $[1,\infty)$ and nested in the sense that $q_{2n+1}<q_{2n+2}$ and $\{q_1,...,q_{2n}\}\cap[q_{2n+1}-\frac{1}{2^n},q_{2n+2}+\frac{1}{2^n}]=\emptyset$ for all $n\ge 0$.  We can then follow the proof of Theorem~\ref{pro} until case 2, the first case being carried out exactly in the same way. 

The second case has to be slightly modified. The technical, but not difficult, part is to ensure that the modified $q_i$'s still form a dense subset of $[1,\infty)$ and that the difference of the left and right derivative at the $q_i$'s does not go to $0$. To do that we take $p_5$  sufficiently close to $p_2$ (in the construction of $\tau_{n+1}$) to have : 
\begin{itemize}
\item $|q_i-\tilde q_i|<\frac{1}{2^n}\inf_{1\le j<j'\le 2n+2}|q_j-q_{j'}|$
\item $|\delta_{n}(q_i)-\delta_{n+1}(\tilde q_i)|<\frac{1}{2^n}\delta_{n}(q_i)$,
\end{itemize}
where $\delta_{n}(q_i)$ stands for the difference between the right and left derivative at $q_i$ of  $\sup_{1\le k\le n}\tau_k$.
The proof of Theorem \ref{3} is then completed in the same way as above. 
\end{Proof}

\bibliographystyle{alpha}
\bibliography{biblio}
\end{document}